\documentclass[letterpaper, 10 pt, conference]{ieeeconf}  

\IEEEoverridecommandlockouts                              
\def\ARXIV{1}

\usepackage[style=ieee,backend=bibtex]{biblatex}
\usepackage{amsmath}
\usepackage{amssymb}
\usepackage[scr=esstix]{mathalpha}
\usepackage{bm}
\usepackage{bbm} 
\usepackage{hyperref}
\hypersetup{colorlinks,
	linkcolor=blue,
	citecolor=blue,
	urlcolor=blue,
	linktocpage,
	plainpages=false}
\usepackage{cancel}
\usepackage{graphicx}
\usepackage{mathtools}
\newtheorem{theorem}{Theorem}

\newtheorem{remark}{Remark}

\newcommand{\R}{\mathbb{R}}

\newcommand{\thetah}{{\hat{\theta}}}
\newcommand{\Thetah}{{\widehat{\Theta}}}

\newcommand{\Ac}{\mathcal{A}}
\newcommand{\Bc}{\mathcal{B}}

\newcommand{\Tc}{\mathcal{T}}

\newcommand{\ts}{\mathscr{dt}} 
\newcommand{\As}{\mathscr{A}}

\newcommand{\Gs}{\mathscr{G}}

\newcommand{\Is}{\mathscr{I}}

\newcommand{\Ls}{\mathscr{L}}

\newcommand{\Ss}{\mathscr{S}}

\newcommand{\Vs}{\mathscr{V}}

\newcommand{\Thetaw}{\widetilde{\Theta}}

\newcommand{\dt}{\text{dt}} 

\newcommand{\tr}{\text{tr}}

\begin{document}

\title{\LARGE \bf
  Sufficient Conditions for Persistency of Excitation with Step and ReLU Activation Functions
}

\author{Tyler Lekang and Andrew Lamperski
  \thanks{This work was supported in part by NSF CMMI-2122856}
      \thanks{T. Lekang and A. Lamperski are with the department of Electrical and
        Computer Engineering, University of Minnesota, Minneapolis,
        MN 55455, USA 
        {\tt\small lekang@umn.edu, alampers@umn.edu}}
}

\maketitle
\thispagestyle{empty}
\pagestyle{empty}

\begin{abstract}
This paper defines geometric criteria which are then used to establish sufficient conditions for
persistency of excitation with vector functions constructed from single hidden-layer neural networks with step or ReLU activation functions.
We show that these conditions hold when employing reference system tracking, as is commonly done
in adaptive control. We demonstrate the results numerically on a system with linearly
parameterized activations of this type and show that the parameter estimates converge to the
true values with the sufficient conditions met.
\end{abstract}

\section{Introduction}
Persistency of excitation is a fundamental concept employed within contexts and applications
related to parameter learning, such as system identification and adaptive control. It is often
discussed, or at least mentioned, in adaptive control textbooks, such as in
\cite{sastry1989adaptive,slotine1991applied,lavretsky2013robust}. \\\\
It was proven in \cite{morgan1977uniform,anderson1977exponential} that persistency of excitation
is necessary and sufficient for the global uniform asymptotic stability of the linear
time-varying (LTV) system
\begin{equation}\label{LTVPE}
\dot{\tilde{\theta}}_t = -V_t\,V_t^\top \tilde{\theta}_t
\end{equation}
where $\tilde{\theta}_t\in\R^n$ is the system state, and $V_t\in\R^{n\times d}$ is a vector
($d=1$) or matrix ($d\geq2$) function of time that is \textit{regulated} (one-sided limits exist
for all $t\in[0,\infty)$).
Consider if $\tilde{\theta}_t=\thetah_t-\theta$ represents the error of a parameter estimate
$\thetah_t$ from some fixed, unknown parameter values $\theta$. Then if $V_t$ is persistently
exciting, the state of this system (ie, the error of the parameter estimates) converges globally
uniformly asymptotically to zero. \\\\
In (non)linear systems with linear parameterizations, the parameter estimate error dynamics
commonly have the form \eqref{LTVPE}. For example, in \cite{slotine1991applied} (sec 8.7), we see
examples of systems that can be formed into a model $y_t = V_t\,\theta$, where the vector $y_t$
and the matrix $V_t$ are measurable, and then using a simple gradient-based update rule for the
parameter estimate $\thetah_t$, within an estimator system $\hat y_t = V_t\,\thetah_t$, gives
exactly these dynamics for the parameter estimate error $\thetah_t-\theta$. Thus, if $V_t$ is
persistently exciting, the parameter estimate error will converge to zero. \\\\
There have been many works since \cite{morgan1977uniform,anderson1977exponential} which utilize
an assumption of persistency of excitation in order to achieve results in parameter learning. In \cite{adetola2008finite}, a
simple parameter learning scheme can be employed for a general class of nonlinear systems which
have some kind of working (nonlinear) feedback controller. An integral condition similar to
persistency of excitation is assumed to be satisfiable. This inspired the work in
\cite{roy2017combined}, which assumes a similar integral condition in order to
identify, and then provide MRAC control for, an unknown MIMO LTI system. And in
\cite{chowdhary2013concurrent,lavretsky2013robust} we see additional examples of MRAC control
which assume persistency of excitation in order to achieve parameter convergence, while in
\cite{annaswamy2021online}, persistently exciting assumptions are made in reinforcement learning
applications. In \cite{willems2005note}, a sufficient condition for windows of observed behavior
of an LTI system (in discrete time) to span the space of possible windows, is for a component
signal (like the input) to be persistently exciting. In \cite{nar2019persistency}, conditions
for neural networks excitation are given to guarantee bounds on the function estimate error. Lastly, in
\cite{panteley2001relaxed} it is proven that for a general class of nonlinear systems which are
feedback linearizable (see \cite{marino1993global,khalil2017high}), global uniform asymptotic
stability can be achieved for linearly paramatrized vector functions meeting relaxed persistency
of excitation conditions. \\\\
However, in these and other works which utilize persistency of excitation assumptions, there is
often no explicit sufficient conditions provided for how to ensure that persistency of
excitation is satisfied. \\\\
On the other hand, there have been some works which do provide these sufficient conditions. A
classic result is that the state of
an LTI system satisfies persistence of excitation if the (stationary) input to that system
contains sufficient frequency content (``sufficient richness", see \cite{boyd1986necessary}). In
\cite{mareels1988persistency,kreisselmeier1990richness} this is extended to certain linear
time-varying systems, while in \cite{lin1999nonlinearities} frequency arguments for sufficient
conditions for excitation are then extended to nonlinear systems in parametric-strict-feedback
form and in \cite{karg2022excitation} to the context of adaptive dynamic programming, in which
optimal control value functions are approximated using polynomial basis functions. In
\cite{padoan2017geometric}, a rank condition is proven sufficient and necessary for the state of
time-invariant systems to be persistently exciting, and in \cite{rueda2021strong}, strong
Lyapunov functions are provided that are equivalent to the persistency of excitation condition.
Finally, and closest in spirit to this paper, in \cite{kurdila1995persistency} a sufficient
condition, based on geometric criteria, is given for satisfying persistency of excitation with
vector functions composed of radial basis functions.\\\\
Our primary contributions are sufficient conditions, based on geometric criteria,
for satisfying persistency of excitation with vector functions $\phi:\R^n\to\R^N$ which are
composed of ReLU or step activation functions, together with  affine transformations of the state
space. We then demonstrate this using a simulated MRAC control application with the parameter
estimates converging to the true values.\\\\
The organization of the remaining parts of the paper are as follows. Section~\ref{sec:notation}
provides preliminary notation. Section~\ref{sec:setup} sets up problem definitions and geometry.
Section~\ref{sec:theory} presents the main theoretical results, while
Section~\ref{sec:numerical} presents numerical results, and closing remarks are given in
Section~\ref{sec:conclusion}.

\section{Notation}
\label{sec:notation}

We interpret $w,x\in\R^n$ as column vectors and denote their inner product as $w^\top\!x$.
Similarly, we denote the product of matrix $W\in\R^{n\times N}$ with vector $x$ as $W^\top\!x$,
which is a length $N$ vector where the $i$th (row) element is the inner product $W_i^\top\!x$.
Index subscripts on vectors and matrices denote the row index, for example $W_i^\top$ is the
$i$th row of $W^\top$. The standard Euclidean norm and Frobenius matrix norm are respectively
denoted $\|w\|_2$ and $\|W\|_F$. The integer set $\{1,\dots,k\}$ is denoted by $[k]$. We use
$t$ subscripts on time-dependent variables to reduce parentheses, for example $\phi(x(t))$ is
instead denoted $\phi(x_t)$. The $i$th row of a time-varying vector or matrix is thus denoted
with $t,i$ subscript. For square $n\times n$ matrix $A$, we use $aI_n\preceq A\preceq bI_n$ to
denote that $\text{eig}(A)\in[a,b]$, where $I_n$ denotes the $n\times n$ identity matrix. We
denote the $n$-length zeros vector as $0_n$.


\section{Setup}
\label{sec:setup}

This section formally defines the nonlinear vector function $\phi$ and what persistency of
excitation means with regards to this definition, then describes the geometry induced on the
state space $\R^n$ by this construction.

\subsection{Nonlinear, Positive Semidefinite Activation Functions}
Let $\phi:\R^n\to\R^N$ be a vector function defined as
\begin{equation}\label{phi}
\phi(x) =
\begin{bmatrix}
\phi_1(x) \\
\vdots \\
\phi_N(x)
\end{bmatrix}
= 
\begin{bmatrix}
\sigma(w_1^\top\!x + b_1) \\
\vdots \\
\sigma(w_N^\top x + b_N)
\end{bmatrix} \ ,
\end{equation}
where $\phi_1,\dots,\phi_N:\R^n\to\R$ are \mbox{composed of \textit{nonlinear}}, piecewise
continuous functions $\sigma:\R\to\R$ together with affine transformations
$w_1^\top\!x + b_1,\dots,w_N^\top x + b_N:\R^n\to\R$. We allow $w_i\in\R^n\setminus\{0_n\}$ and
$b_i\in\R$ to be arbitrary for all $i\in[N]$, except we assume each $w_i^\top\!x + b_i = 0$
hyperplane in $\R^n$ is unique \mbox{with dimension $n-1$.} Let $W = \begin{bmatrix}w_1 &\cdots
& w_N\end{bmatrix}$ and $b=\begin{bmatrix}b_1 & \cdots & b_N \end{bmatrix}^\top $. For any
$\Ss\subset [N]$, let $W^\top_{\Ss}$ and be the submatrix of $W^\top$ with rows given by
$w_i^\top$ for $i\in\Ss$. Define $b_{\Ss}$ similarly for $b$. \\\\
Note then that \eqref{phi} is equivalent to the output of a single hidden neural network layer
with $N$ neurons fully connected to the input $x\in\R^n$, having nonlinear (eg, ReLU)
activations, and being initialized with weights and biases defining unique hyperplanes. Hence
why we refer to $\sigma$ as an \textit{activation} function or simply an activation.
This paper will focus on the following activations:
\begin{align}
\sigma_{cs}(y) &= 
\begin{cases*}\label{step}
0 \quad\text{if } y\leq0 \\
c \quad\text{if } y>0
\end{cases*}
\qquad \text{(scaled step)} \\
\sigma_r(y) &= 
\begin{cases*}\label{ReLU}
0 \quad\text{if } y\leq0 \\
y \quad\text{if } y>0
\end{cases*}
\qquad \text{(ReLU)} \ \ ,
\end{align}
where $c>0$ is an arbitrary positive scalar.

\subsection{Persistency of Excitation}
Let $x:[0,\infty)\to\R^n$ be some continuous trajectory in the state space. Then, $\phi(x_t)$ is
a piecewise continuous (and regulated) vector function of time, mapping $[0,\infty)\to\R^N$. For
any time window $t\in[\tau,\tau+T]$ with $\tau\geq0$ and $T>0$, the integral
\begin{equation}\label{gram}
\int_{\tau}^{\tau+T} \phi(x_t)\,\phi(x_t)^\top\,\dt
\end{equation}
defines a $N\times N$ \textit{Gramian matrix} since the corresponding \mbox{$i,j$-th} entry
$\int_{\tau}^{\tau+T} \phi_i(x_t)\,\phi_j(x_t)\,\dt$ is an inner product of the
composition functions $\phi_1(x_t),\dots,\phi_N(x_t):[0,\infty)\to\R$. Gramian matrices are
always positive semidefinite, which can be shown using the bi-linearity of inner products.\\\\
\textit{Persistency of excitation} is the requirement that the Gramian matrix \eqref{gram} 
must be strictly positive definite, with eigenvalues in some bounded interval
$[\alpha_1,\alpha_2]$, over all shifts $\tau\geq0$ of the sliding time window $[\tau,\tau+T]$
for some window length $T>0$. Formally, persistency of excitation requires existence of
constants $\alpha_1,\alpha_2,T>0$ such that
\begin{equation}\label{PE}
\alpha_1 I_N\ \preceq \ \int_{\tau}^{\tau+T} \phi(x_t)\,\phi(x_t)^\top\,\dt \ \preceq \ \alpha_2 I_N
\end{equation}
holds for all $\tau\geq0$. An equivalent scalar requirement is that
\begin{equation}\label{PEvec}
\alpha_1\|v\|_2^2 \ \leq \ 
\int_{\tau}^{\tau+T} (v^\top \phi(x_t))^2\,\dt \ \leq \
\alpha_2\|v\|_2^2
\end{equation}
must hold for all $v\in\R^N$ and $\tau\geq0$. This follows since $v^\top
\phi(x_t)=\phi(x_t)^\top v$.\\\\
Note that if these hold for some $T>0$, then they also hold for all $\widetilde{T}>T$. This
follows because the integrals can be broken into a sum of two integrals over
$t\in[\tau,\tau+T]$ and $t\in[\tau+T,\tau+\widetilde{T}]$, with the former being strictly
positive (definite) and the latter positive semidefinite (nonegative).

\subsection{Activation Geometry}
We define the following subsets for all activations $i\in[N]$:
\begin{align*}
X_i^+ &= \{x\in\R^n \ | \ w_i^\top\!x + b_i > 0 \} \qquad \text{(active)} \\
X_i^\circ &= \{x\in\R^n \ | \ w_i^\top\!x + b_i \leq 0 \} \qquad \text{(zero)} \ \ .
\end{align*}
Each is a half-space of $\R^n$ formed by the hyperplane $w_i^\top\!x + b_i = 0$. These are then
used to define the \textit{activation regions} $\Ac_j$, with indices corresponding to a binary
string indicating which active (binary 1) and zero (binary 0) half-spaces are in the
intersection:
\begin{align*}
\Ac_0 &= X_N^\circ\ \cap\ X_{N-1}^\circ\ \cap\ \dots\ \cap\ X_2^\circ \cap\ X_1^\circ \\
\Ac_1 &= X_N^\circ\ \cap\ X_{N-1}^\circ\ \cap\ \dots\ \cap\ X_2^\circ \cap\ X_1^+ \\
&\hspace{6pt}\vdots \\
\Ac_{2^N-1} &= X_N^+\ \cap\ X_{N-1}^+\ \cap\ \dots\ \cap\ X_2^+ \cap\ X_1^+ \ \ .
\end{align*}
These partition $\R^n$ into, at most, $2^N$ convex polytopes. It is likely that some of the
regions will be infeasible ($\Ac_j=\emptyset$). For example, in the $n=1$ case, there are always
only $N+1$ feasible activation regions, since there will be $N$ unique (by assumption) points
partitioning the $\R$ line. In higher dimensions, more feasible regions are possible.\\\\
For all $j\in\{0,\dots,2N-1\}$, we define the \textit{active set} $\Ss_j$ as
\begin{equation*}
\Ss_j = \{i\in[N] \ | \ w_i^\top\!x + b_i > 0 \quad \forall x\in\Ac_j \} \ \ .
\end{equation*}
This captures which of the $N$ activations are active in a particular activation region.
For any two activation regions $j,k\in\{0,\dots,2N-1\}$ with a nonempty intersection
$\Ac_j\cap\Ac_k\neq\emptyset$, their intersection defines a \textit{border} between the two
regions. We define a \textit{nondegenerate border} to mean $\dim(\Ac_j \cap \Ac_k)=n-1$. In this
case, only one activation is different (active to zero or vice-versa) between those regions.
This is because the borders between activation regions must be (subsets of) the hyperplanes that
define the half-spaces $X_i$, and the intersection of more than one unique hyperplane with
dimension $n-1$ must have a dimension \textit{less than} $n-1$. Thus, a nondegenerate border
must be a subset of (or equal to) a single hyperplane, meaning only one $X_i$ half-space can
flip from active to zero or vice-versa. We then define a \textit{degenerate border} to mean
$\dim(\Ac_j \cap \Ac_k)<n-1$. In this case, the border is a (subset of) the lower dimensional
intersection of two or more unique hyperplanes. Thus, multiple activations are different between
the regions.\\\\
Now consider a continuous state trajectory $x_t\in\R^n$ visiting a sequence of activation
regions over the time window $[\tau,\tau+T]$, for some $\tau\geq0$ and $T>0$. Assume $x_t$ only
crosses nondegenerate borders and that the number of regions visited is $\Ls\geq2$. Let us denote
$\As_1,\dots,\As_\Ls$ as the activation region indices of the visited sequence and
$\Is_1,\dots,\Is_{\Ls-1}$ as the sequence of activation indices corresponding to the
hyperplanes crossed in order to visit that sequence of activation regions. For all $s\in[\Ls]$,
we define time window subsets $\Tc_s \subset [\tau,\tau+T]$ as
\begin{equation}\label{eq:timePartition}
\Tc_s = \{t\in[\tau,\tau+T] \ | \ x_t \in \Ac_{\As_s} \} \ \ .
\end{equation}
Since we assume $x_t$ crosses only nondegenerate borders when visiting the activation
regions during the time window, we have by definition over all $s\in[\Ls-1]$ that
\begin{equation}\label{activeSetInduction}
\Ss_{\As_{s+1}} =
\begin{cases}
\Ss_{\As_s} \bigcup \ \{\Is_s\} \qquad \text{if}\ \ \Ac_{\As_{s+1}}\subseteq
X_{\Is_s}^+ \\
\Ss_{\As_s} \setminus\ \{\Is_s\} \qquad \text{if}\ \ \Ac_{\As_{s+1}}\subseteq
X_{\Is_s}^\circ
\end{cases} .
\end{equation}


\section{Theoretical Results}
\label{sec:theory}

This section presents our main theoretical results, which provide sufficient conditions for
satisfying persistency of excitation with (scaled) step or ReLU activations.

\subsection{Main Results}

\begin{theorem} \label{thm:step} \textit{
Let state trajectory $x_t$ be continuous and stay within some compact set $\Bc\subset\R^n$ for
all $t\geq0$, and let $\phi(x_t)=[\sigma_{cs}(w_1^\top x_t+b_1),\dots,\sigma_{cs}(w_N^\top
x_t+b_N)]$ be composed of (scaled) step functions \eqref{step} with positive scalars
$c=[c_1,\dots,c_N]^\top$ together with $N$ unique affine transformations of $\R^n$ according to
$w_1,\dots,w_N\in\R^n\setminus\{0_n\}$ and $b_1,\dots,b_N\in\R$. If $x_t$ over $t\geq0$ is such
that there exists a window length $T^*>0$ whereby the sequence $s\in[\Ls]$ of activation regions
visited during any shift $\tau\geq0$ of the time window $[\tau,\tau+T^*]$ always satisfies the
following two conditions:\\
\phantom{ }\hspace{2pt} 1. all $i\in[N]$ hyperplanes $w_i^\top x_t+b_i=0$ are crossed \\
\phantom{ }\hspace{2pt} 2. only nondegenerate borders are crossed, \\
then $\phi(x_t)$ satisfies the persistency of excitation conditions \eqref{PE} and \eqref{PEvec}.
}
\end{theorem}
\begin{proof}
\if\ARXIV1
Given in Appendix~\ref{appThmProofs}.
\fi
\if\ARXIV0
Given in Appendix I of \cite{lekang2022sufficientFULL}.
\fi
\end{proof}
\vphantom{~}
\begin{theorem} \label{thm:ReLU} \textit{
Let state trajectory $x_t$ be continuous and stay within some compact set $\Bc\subset\R^n$ for
all $t\geq0$, and let $\phi(x_t)=[\sigma_r(w_1^\top x_t+b_1),\dots,\sigma_r(w_N^\top x_t+b_N)]$
be comprised of ReLU functions \eqref{ReLU} together with $N$ unique affine transformations of
$\R^n$ according to $w_1,\dots,w_N\in\R^n\setminus\{0_n\}$ and $b_1,\dots,b_N\in\R$. If $x_t$
over $t\geq0$ is such that there exists a window length $T^*>0$ whereby the sequence $s\in[\Ls]$
of activation regions visited during any shift $\tau\geq0$ of the time window
$[\tau,\tau+T^*]$ always satisfies the following three conditions: \\
\phantom{ }\hspace{2pt} 1. all $i\in[N]$ hyperplanes $w_i^\top x_t+b_i=0$ are crossed \\
\phantom{ }\hspace{2pt} 2. only nondegenerate borders are crossed \\
\phantom{ }\hspace{2pt} 3. for each $s\in[\Ls]$, there are times $t_1,\hat t_1,\ldots,t_M,
\hat t_M\in \Tc_s$ such that
$$
\mathrm{rank}(W^\top_{\Ac_{\As_s}})=
\mathrm{rank}\left( W^\top_{\Ac_{\As_s}}
[x_{t_1}-x_{\hat t_1} \ \cdots \ x_{t_M}-x_{\hat t_M}]\right),
$$
then $\phi(x_t)$ satisfies the persistency of excitation conditions \eqref{PE} and \eqref{PEvec}.
}
\end{theorem}
\begin{proof}
\if\ARXIV1
Given in Appendix~\ref{appThmProofs}.
\fi
\if\ARXIV0
Given in Appendix I of \cite{lekang2022sufficientFULL}.
\fi
\end{proof}
Note that a sufficient condition for property (iii) is that
$\begin{bmatrix}x_{t_1}-x_{\hat t_1} & \cdots & x_{t_M}-x_{\hat t_M}\end{bmatrix}$ has rank $n$.
This can be achieved, for example, if $x_t$ is the state trajectory of a system that satisfies a
suitable (local) controllability property.

\subsection{Proof Sketch}
Both proofs rely on the same overall contradiction method. That is, we assume there exists a
\underline{nonzero} vector $v\in\R^N\setminus\{0_N\}$ such that
\begin{equation}\label{nonPEvec0}
v^\top\left( \int_{\tau}^{\tau+T} \phi(x_t)\,\phi(x_t)^\top\dt\right) v \ = \ 0 \ \ .
\end{equation}
We show that if the state trajectory $x_t$ meets certain requirements over
any shift $\tau\geq0$ of the window $[\tau,\tau+T^*]$, for some window length $T^*>0$, then in
fact \eqref{nonPEvec0} can only hold if $v=0_N$. This is a contradiction and thus proves that
the LHS integral must be strictly positive definite, uniformly over all windows
$[\tau,\tau+T^*]$ for all $\tau\geq0$.

\if\ARXIV1
\subsection{Extension to Other Activation Functions}

In Appendix~\ref{appThmProofs} we discuss how these methods could be applied similarly to other
nonlinear, positive semidefinite activation functions $\sigma$, such that $\sigma(y)=0$ for all
$y\leq 0$ and $\sigma(y) >0$ for all $y>0$. 
\fi


\section{Numerical Results}
\label{sec:numerical}

In this section we provide a numerical simulation\footnote{All code is available at:
\href{https://github.com/tylerlekang/CDC2022}{https://github.com/tylerlekang/CDC2022}} of the
theoretical results, using a MRAC application which is a variation on the setup from Chapter 9
of \cite{lavretsky2013robust}.
The plant and reference systems have $n=2$ states, allowing convenient visualization of the
hyperplanes and state space.

\subsection{Setup}

The plant is given by
\begin{equation}\label{plant}
\dot{x}_t = A\,x_t + B\big(u_t + \Theta^\top\!\phi(x_t)\big) \ \ ,
\end{equation}
where $A$ is a known $n\times n$ state matrix for the plant state $x_t\in\R^n$, $B$ is a known
$n\times\ell$ input matrix for the input $u_t\in\R^\ell$, and $\Theta$ is an \underline{unknown}
$N\times\ell$ matrix which linearly parameterizes the known vector function $\phi:\R^n\to\R^N$
defined by \eqref{phi}. The setup in \cite{lavretsky2013robust} also includes an unknown
diagonal scaling matrix $\Lambda$, such that the overall input matrix is $B\Lambda$. We have
omitted this for simplicity.\\\\
The control input $u_t$ will be designed in order to force the plant states $x_t$ to track the
states of a reference system $x_t^r$ that is driven by a bounded reference input $r_t$. The
reference system is given by
\begin{equation}
\dot{x}_t^r = A_r\,x_t^r + B_r\,r_t \ \ ,
\end{equation}
where $A_r$ and $B_r$ are known reference matrices, with $A_r$ Hurwitz, and $r_t\in\R^\ell$ is a
bounded reference input.\\\\
We assume there exists an $n\times\ell$ matrix of feedback gains $K_x$ and an $\ell\times\ell$
matrix of feedforward gains $K_r$ satisfying the \textit{matching conditions}
\begin{equation}\label{matching}
A + BK_x^\top = A_r \quad\text{and}\quad BK_r^\top = B_r \ \ .
\end{equation}
The setup in \cite{lavretsky2013robust} has $A$ and $\Lambda$ as unknown, and thus $K_x$ and
$K_r$ need to be estimated. For this simulation, we will assume that $K_x$ and $K_r$ can be
directly calculated from known $A$ and $B$, and used within the control law.\\\\
Next, we introduce parameter estimates $\Thetah_t$, which will be dynamically updated to
estimate true parameter values $\Theta$. Thus, by applying to the plant \eqref{plantSys}
the feedback control law $u_t=K_x^\top x_t+K_r^\top r_t - \Thetah_t^\top\!\phi(x_t)$, the plant
dynamics become
\begin{equation}
\dot{x}_t = A_r\,x_t + B_r\,r_t - B\,(\Thetah_t-\Theta)^\top\!\phi(x_t) \ \ .
\end{equation}
This in turn gives the dynamics of the state tracking error $e_t = x_t-x_t^r$ as
\begin{equation}\label{stTrkErrDyn}
\dot{e}_t = \dot{x}_t-\dot{x}_t^r = A_r\,e_t - B\,(\Thetah_t-\Theta)^\top\!\phi(x_t) \ \ .
\end{equation}
In \cite{lavretsky2013robust}, it is then shown that these state tracking error dynamics
$\dot{e}_t$ are globally uniformly asymptotically stable, such that
$\lim_{t\to\infty}\|e_t\|_2=0$, if the parameter estimates are dynamically updated as
\begin{equation}\label{update}
\dot{\Thetah}_t = \Gamma\,\phi(x_t)\,e_t^\top P_xB \ \ .
\end{equation}
This is shown by analyzing the Lyapunov function $\Vs =
e_t^\top P_xe_t + \tr\big((\Thetah_t-\Theta)^\top\Gamma^{-1}(\Thetah_t-\Theta)\big)$, along with
using Barbalat's lemma, such that the update rule \eqref{update} results in
$\dot{\Vs}\leq-e_t^\top Q_xe_t$ for all values of $e_t$ and $\Thetah_t-\Theta$. Here, $P_x$ is
the unique symmetric, positive-definite $n\times n$ matrix that solves the algebraic Lyapunov
equation $P_xA_r + A_r^\top P_x=-Q_x$ for some symmetric, positive-definite $n\times n$ matrix
$Q_x$, and \textit{adaptation rates} $\Gamma$ is some symmetric, positive-definite $N\times N$
matrix, where we denote $\|\Gamma\|_F:=\Gs$.

\begin{remark}
In the case where $K_x$ and $K_r$ are also being estimated by $\hat{K}_{x,t}$ and
$\hat{K}_{r,t}$ respectively, the true values can be appended to $\Theta$, the estimates can be
appended to $\Thetah_t$, and an overall vector function $\Phi(x_t,r_t)$ which combines $x_t$,
$r_t$, and $\phi(x_t)$, can be formed. However, such $\Phi$ then do not strictly meet the
definition of \eqref{phi}, and are thus beyond the scope of this paper.
\end{remark}

\subsection{Persistency of Excitation}
The dynamic update rule \eqref{update} only guarantees asymptotic convergence of the state
tracking error $e_t$ to zero. We now show that the parameter estimation error, which we will
denote compactly as $\Thetaw_t = \Thetah_t-\Theta$, also goes to zero if $\phi(x_t)$ is
persistently exciting.\\\\
Since $A_r$ is Hurwitz, it is guaranteed invertible. And so, from \eqref{stTrkErrDyn} we
have
\begin{equation*}
e_t^\top P_xB = \dot{e}_t^\top {A_r^{-1}}^\top P_xB +
\phi(x_t)^\top \Thetaw_t\,B^\top {A_r^{-1}}^\top P_xB
\end{equation*}
and then combined with \eqref{update} we get
\begin{align*}
&\dot\Thetaw_t = \dot{\Thetah}_t = \Gamma\,\phi(x_t)\,e_t^\top P_xB = \\
& \Gamma\,\phi(x_t)\phi(x_t)^\top \Thetaw_t\,B^\top {A_r^{-1}}^\top P_xB +
\Gamma\,\phi(x_t)\,\dot{e}_t^\top {A_r^{-1}}^\top P_xB \ .
\end{align*}
Note that the second term asymptotically goes to zero with $\dot{e}_t$ and we have in the first
term $B^\top{A_r^{-1}}^\top P_xB=$ $-\frac{1}{2}B^\top {A_r^{-1}}^\top Q_x\, A_r^{-1}B$ 
by the Lyapunov equation definition for $P_x$.\\\\
Let us now restrict to the case of $\ell=1$. Since $Q_x$ is positive-definite we have the
positive scalar $c=\frac{1}{2}B^\top {A_r^{-1}}^\top Q_x\, A_r^{-1}B$, and obtain that the
dynamics of the parameter estimation error vector $\Thetaw_t\in\R^N$ asymptotically tend to the
linear time-varying system
\begin{equation}\label{paramEstDyn}
\dot\Thetaw_t = -\Gamma_c\,\phi(x_t)\phi(x_t)^\top \Thetaw_t \ \ .
\end{equation}
$\Gamma_c = c\,\Gamma$ is symmetric, positive-definite with $\|\Gamma_c\|_F=c\Gs$.\\\\
Therefore, if $\phi(x_t)$ satisfies the persistency of excitation condition \eqref{PE},
then the parameter estimation error dynamics \eqref{paramEstDyn} are globally uniformly
asymptotically stable such that we have $\lim_{t\to\infty}\|\Thetaw_t\|_2=0$. We show explicitly
in
\if\ARXIV1
Appendix~\ref{appLTVanderson}
\fi
\if\ARXIV0
Appendix II of \cite{lekang2022sufficientFULL}
\fi
how this follows from Lemmas 1 and 2 and Theorem 1 of
\cite{anderson1977exponential}.

\subsection{Simulation}

We perform the numerical simulation in the $x_t\in\R^2$ state space.
Euler integration with a timestep of $\ts=0.001$ sec. is used to obtain all solutions to
differential equations.
We use a basic controllable canonical form for the plant:
\begin{equation}\label{plantSys}
\dot{x}_t =
\begin{bmatrix}
0 & 1 \\
-a_1 & 2a_2
\end{bmatrix}\!
\begin{bmatrix}
x_{t,1}\\
x_{t,2}
\end{bmatrix} +
\begin{bmatrix}
0\\
\beta
\end{bmatrix}
\big(u_t + \Theta^\top\!\phi(x_t)\big) \ \ ,
\end{equation}
where $a_1,a_2,\beta>0$ are positive scalars, $u_t\in\R$ is the $\ell=1$ dimensional input, and
$\Theta$ is a fixed vector of nonzero scalars which linearly parameterizes the known nonlinear
vector function $\phi(x_t)=[\sigma_R(w_1^\top x_t+b_1)\ \cdots\ \sigma_R(w_4^\top x_t+b_4)]^\top$
with ReLU activation functions,
meeting the definition \eqref{phi}. This plant model is an $n=2$ example of a general
class of single-input systems with dynamics characterized by an $n$th order nonlinear ordinary
differential equation. See section 9.5 in \cite{lavretsky2013robust} and section 4.1 in
\cite{khalil2017high} for examples of physical systems in this form using various nonlinear
functions.\\\\
The reference system is a unity-gain damped harmonic oscillator with a bounded reference input
$r_t\in\R$:
\begin{equation}\label{refSys}
\dot{x}_t^r =
\begin{bmatrix}
0 & 1 \\
-\omega_0^2 & -2\xi\omega_0
\end{bmatrix}\!
\begin{bmatrix}
x_{t,1}^r\\
x_{t,2}^r
\end{bmatrix} +
\begin{bmatrix}
0\\
\omega_0^2
\end{bmatrix}
r_t \ \ ,
\end{equation}
where $\omega_0,\xi>0$ are the natural frequency and damping ratio. This gives the plant and
reference eigenvalues as $\lambda = a_2 \pm \sqrt{a_2^2 - a_1}$ and $\lambda^r = (-\xi \pm
\sqrt{\xi^2 - 1\,}\,)\,\omega_0$, thus $A$ is unstable with oscillations if $a_1>a_2^2$ and
$A_r$ is always Hurwitz, and without oscillations if $\xi\geq1$. Thus, we can always directly
calculate the required feedback and feedforward gains that satisfy the matching conditions
\eqref{matching} as
$$
K_x^\top=\left[-\frac{\omega_0^2-a_1}{\beta}\ \ -\frac{2\xi\omega_0+2a_2}{\beta}\right]
\ \ \text{and} \ \
K_r^\top=\left[\,\frac{\omega_0^2}{\beta}\,\right] \ .
$$\\
For the plant system, we use $a_1=2$, $a_2=0.5$, $\beta=0.75$. This results in an unstable,
oscillatory $A$ matrix with $\text{eig}(A)=0.5 \pm j\sqrt{1.75\,}\,$. The reference system uses
$\omega_0=2$ and $\xi=1$, ensuring that $A_r$ is Hurwitz and nonoscillatory with
$\text{eig}(A_r) = -2\,$. The matching conditions are thus satisfied with $K_x^\top =
[-2.667\ \ -6.667]$ and $K_r^\top=[5.333]$. For the parameter estimate dynamic updates we use:
\begin{equation*}
\Gamma =
\begin{bmatrix}
5 & 0 & 0 & 0 \\
0 & 1 & 0 & 0 \\
0 & 0 & 5 & 0 \\
0 & 0 & 0 & 2
\end{bmatrix}
Q_x =
\begin{bmatrix}
1 & 0 \\
0 & 10
\end{bmatrix}
P_x =
\begin{bmatrix}
5.625 & 0.125 \\
0.125 & 1.281
\end{bmatrix},
\end{equation*}
and for the linearly parameterized nonlinear vector function $\Theta^\top\phi(x_t)$ we use:
\begin{equation*}
\Theta = 
\begin{bmatrix}
-1.2 \\
2.7 \\
0.8 \\
-3.2
\end{bmatrix}
\ \ 
W^\top = 
\begin{bmatrix}
2 & 1 \\
1 & -2 \\
1.5 & -0.5 \\
0.5 & 2
\end{bmatrix}
\ \ 
b = 
\begin{bmatrix}
1 \\
2 \\
2.5 \\
3
\end{bmatrix} \ \ .
\end{equation*}
\if\ARXIV1
In Figure~\ref{fig:pe:R2hyps} we plot points to highlight each activation region, partitioning
the state space $\R^2$.\\\\
\fi
Defining parameter estimates $\Thetah_t=[\thetah_{t,1}\ \cdots\ \thetah_{t,4}]^\top$ and
using the known $K_x$ and $K_r$, we apply the feedback control law $u_t=K^x x_t+K^r r_t -
\Thetah_t^\top\!\phi(x_t)$ to the plant and update the parameter estimates according to the
update law \eqref{update}.
\if\ARXIV0
\newpage\noindent
\fi
We use the following bounded reference inputs to drive the reference system as two different
scenarios:
\begin{equation*}
r_t^{(1)} = 10\sin(0.5t) \qquad r_t^{(2)} = 40 + \sum_{k=1}^2 10\sin(0.25k\,t) \ .
\end{equation*}
The resulting plant and reference state trajectories $x_t$ and $x_t^r$ are plotted in state
space along with the hyperplanes $W^\top x + b=0_4$, in Figure~\ref{fig:phase} and
Figure~\ref{fig:phase2} respectively for the two scenarios.
\if\ARXIV0
We see clearly that the limit cycles in both cases stay within a compact set $\Bc\in\R^2$, and
so the ReLU activations within $\phi$ are bounded on this $\Bc$, and the trajectories never
maintain a linear path.
The parameter estimation error $\|\Thetah_t-\Theta\|_2$ for both cases is plotted over simulation
time in Fig.~\ref{fig:thEstConv}. For the case with $r_t^{(1)}$, we see that the error converges
to zero as the parameter estimators converge to the true values, while for the case with
$r_t^{(2)}$ the error does not converge and is well above zero.
\fi
\if\ARXIV1
We see clearly that the limit cycles in both cases stay within a compact set $\Bc\in\R^2$, and
so the ReLU activations within $\phi$ are bounded on this $\Bc$, and the trajectories never
maintain a linear path.\\\\
In the first case, the state trajectories cross all $N=4$ hyperplanes at nondegenerate borders
persistently, as well as never being parallel to any subset of active hyperplanes, thus meeting
the criteria of Theorem~\ref{thm:ReLU}. In the second case, the limit cycle only crosses two of
the hyperplanes persistently.\\\\
The parameter estimation error $\|\Thetah_t-\Theta\|_2$ for both cases is plotted over simulation
time in Fig.~\ref{fig:thEstConv}. For the case with $r_t^{(1)}$, we see that the error converges
to zero as the parameter estimators converge to the true values, while for the case with
$r_t^{(2)}$ the error does not converge and is well above zero.
\fi

\if\ARXIV1
\begin{figure}[h]
\centering
\includegraphics[width=1.0\columnwidth]{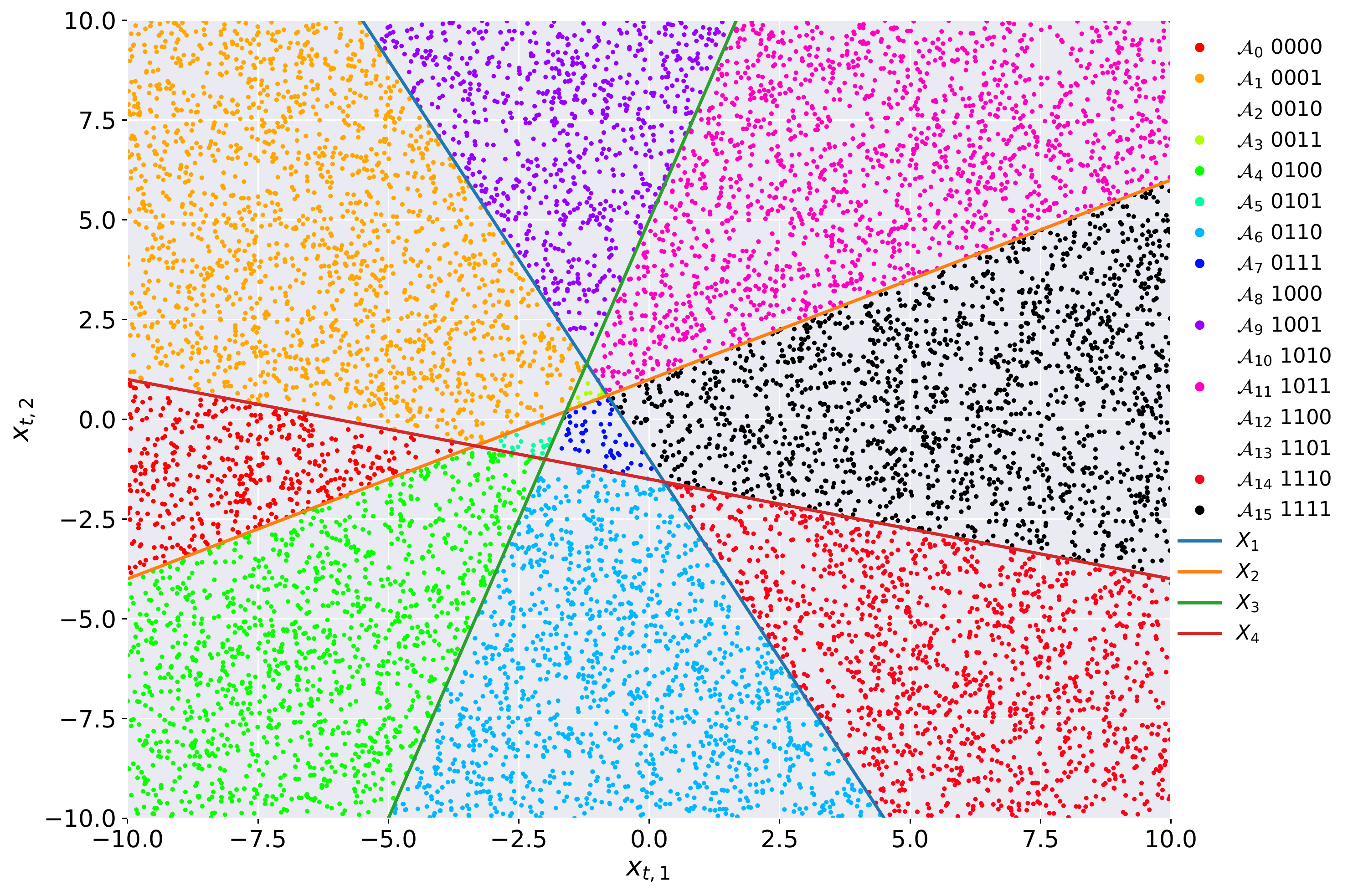}
\caption{\label{fig:pe:R2hyps}State space $x\in\R^2$ showing $N=4$ hyper planes
$W^\top x+b=0_4$. There are $11$ feasible activation regions $\Ac_j$ (color dots) and $5$
infeasible regions (no dots).}
\end{figure}

\begin{figure}[h]
\centering
\includegraphics[width=1.0\columnwidth]{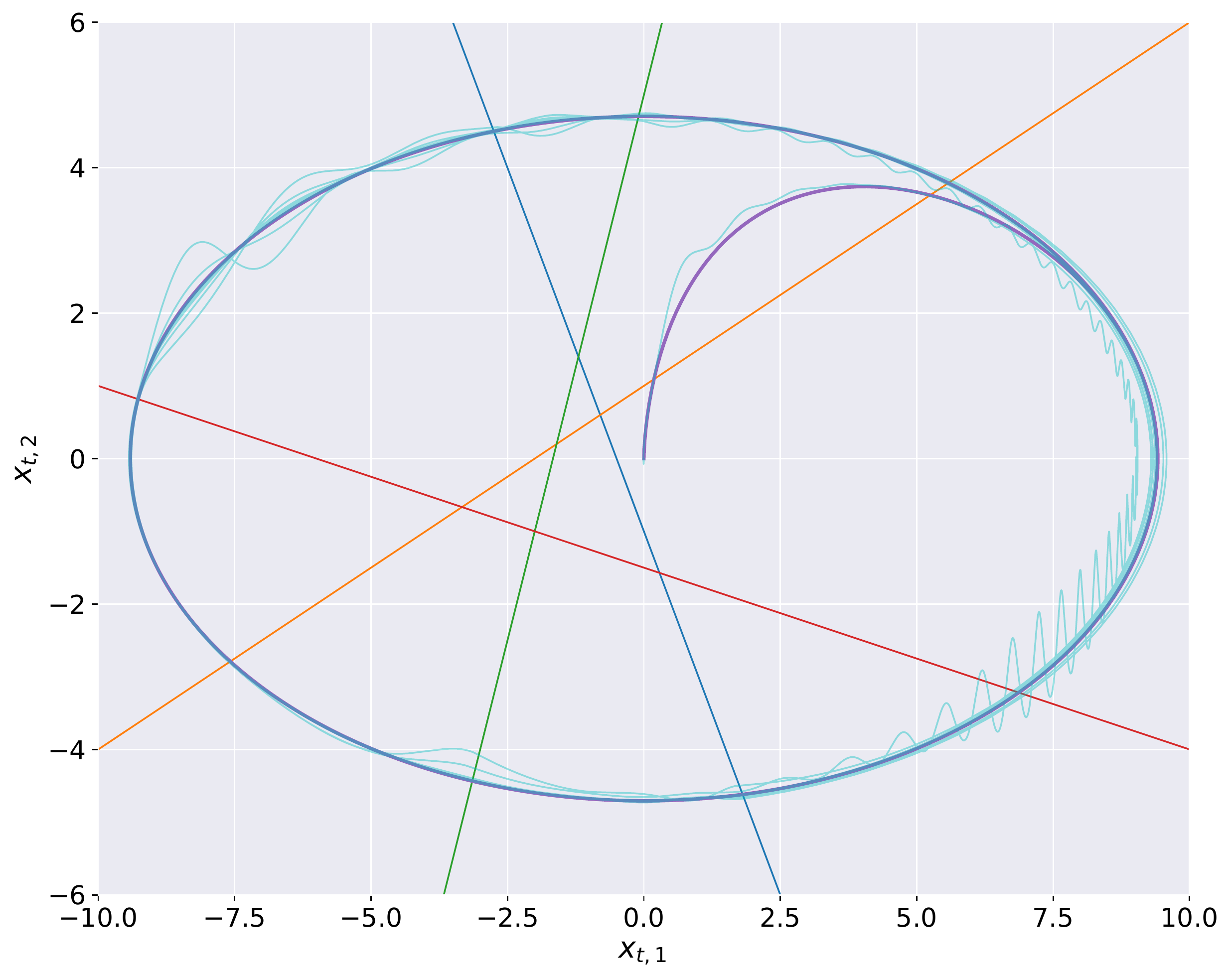}
\caption{\label{fig:phase}Phase plot of $x_t$ (cyan), tracking into $x_t^r$ (magenta) driven by
$r_t^{(1)}$, and crossing all $N$ hyperplanes at nondegenerate borders along the limit
cycle.}
\end{figure}

\begin{figure}[h]
\centering
\includegraphics[width=1.0\columnwidth]{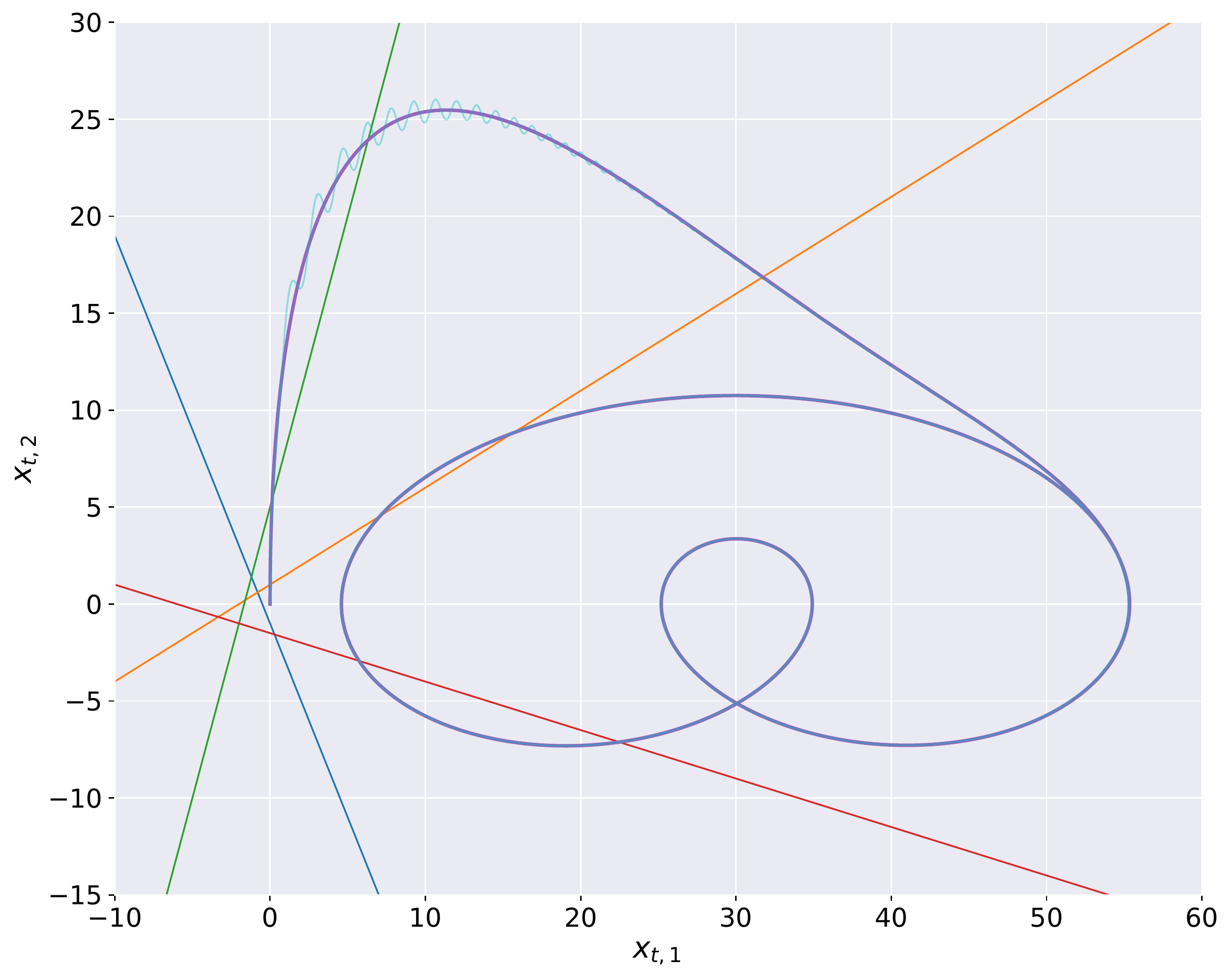}
\caption{\label{fig:phase2}Phase plot of $x_t$ (cyan), tracking into $x_t^r$ (magenta) driven by
$r_t^{(2)}$, but does \textit{not} cross all $N$ hyperplanes along the limit cycle.}
\end{figure}

\begin{figure}[h!]
\centering
\includegraphics[width=1.0\columnwidth]{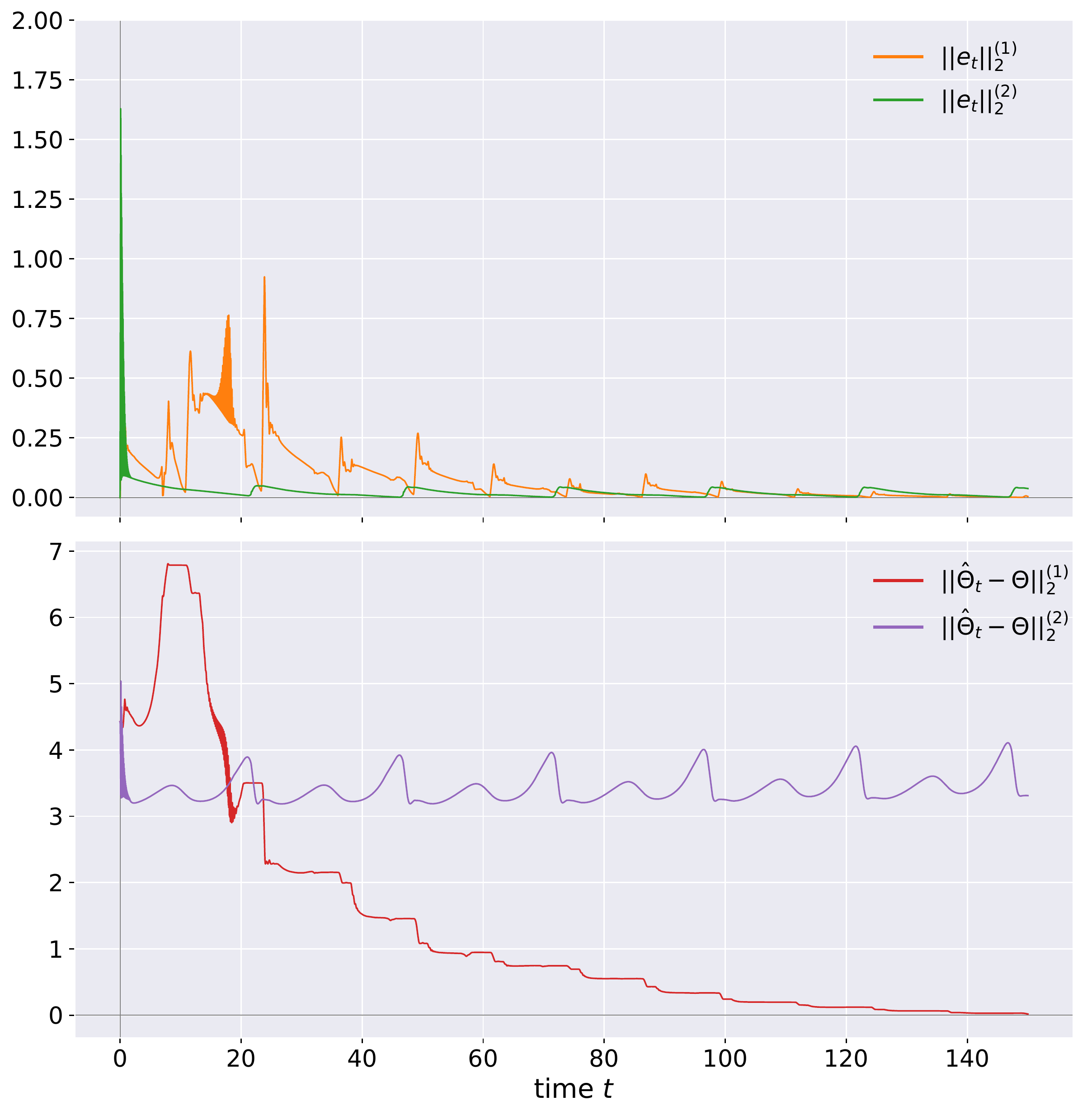}
\caption{\label{fig:thEstConv}Norm of state tracking error $\|e_t\|_2$ converges to zero in both
cases, while the parameter estimation error $\|\Thetah_t-\Theta\|_2$ converges to zero for the
case with $r_t^{(1)}$ and but not for the case with $r_t^{(2)}$.}
\end{figure}
\fi

\if\ARXIV0
\begin{figure}[h]
\centering
\includegraphics[width=0.7\columnwidth]{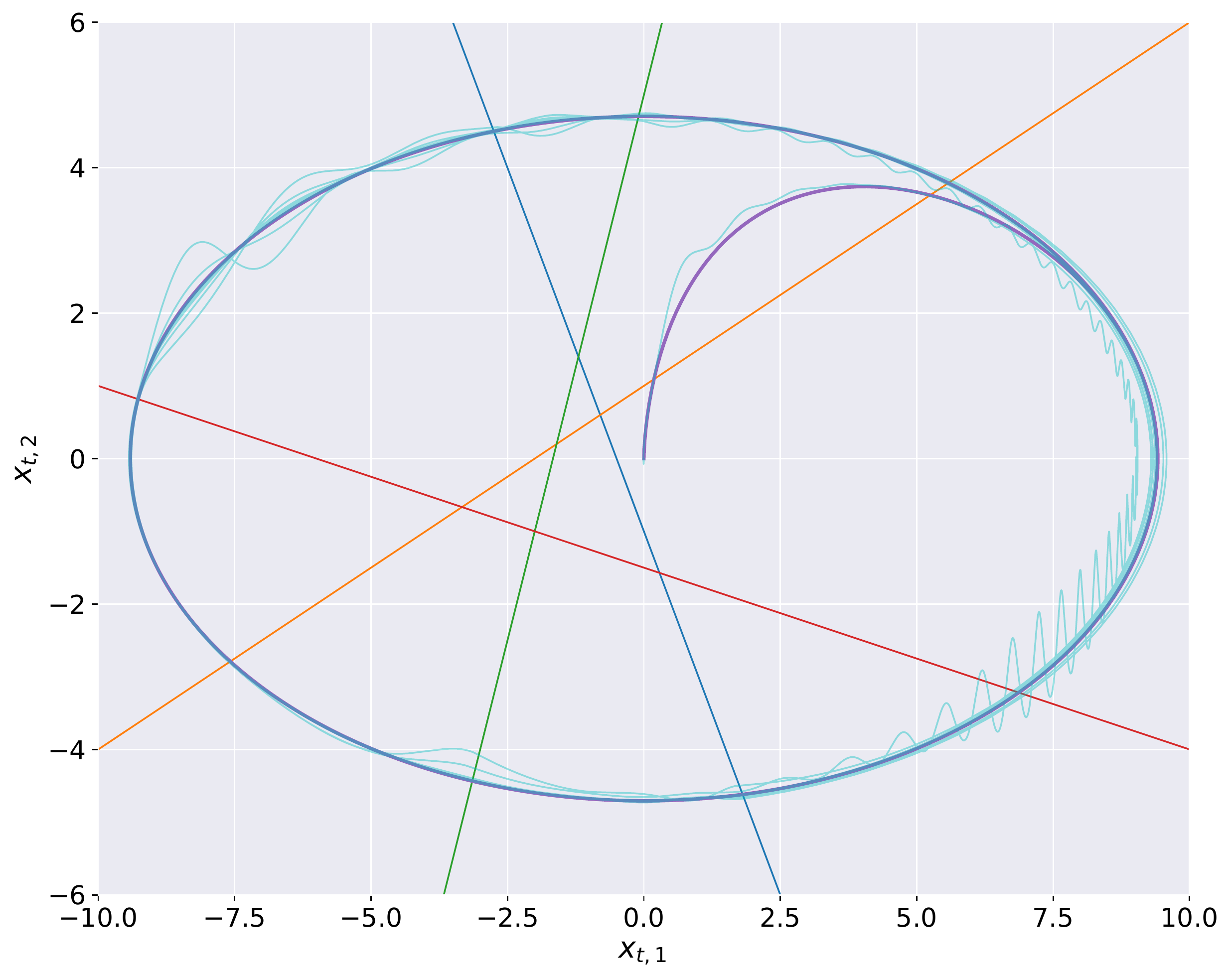}
\caption{\label{fig:phase}Phase plot of $x_t$ (cyan), tracking into $x_t^r$ (magenta) driven by
$r_t^{(1)}$. Crosses all $N$ hyperplanes at nondegenerate borders along the limit cycle.}
\end{figure}
\begin{figure}[h]
\centering
\includegraphics[width=0.7\columnwidth]{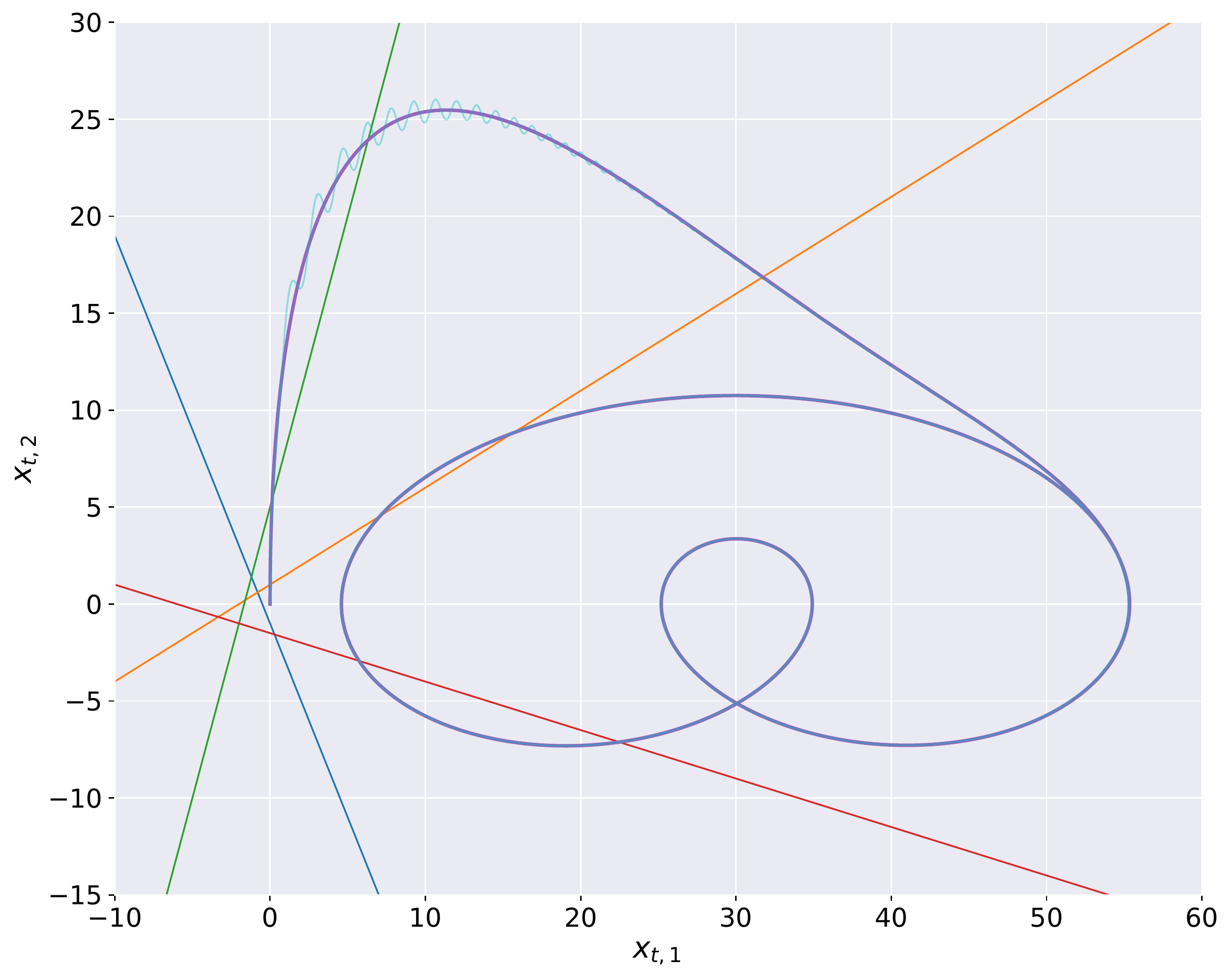}
\caption{\label{fig:phase2}Phase plot of $x_t$ (cyan), tracking into $x_t^r$ (magenta) driven by
$r_t^{(2)}$. Does \textit{not} cross all $N$ hyperplanes along the limit cycle.}
\end{figure}

\begin{figure}[h]
\centering
\includegraphics[width=1.0\columnwidth]{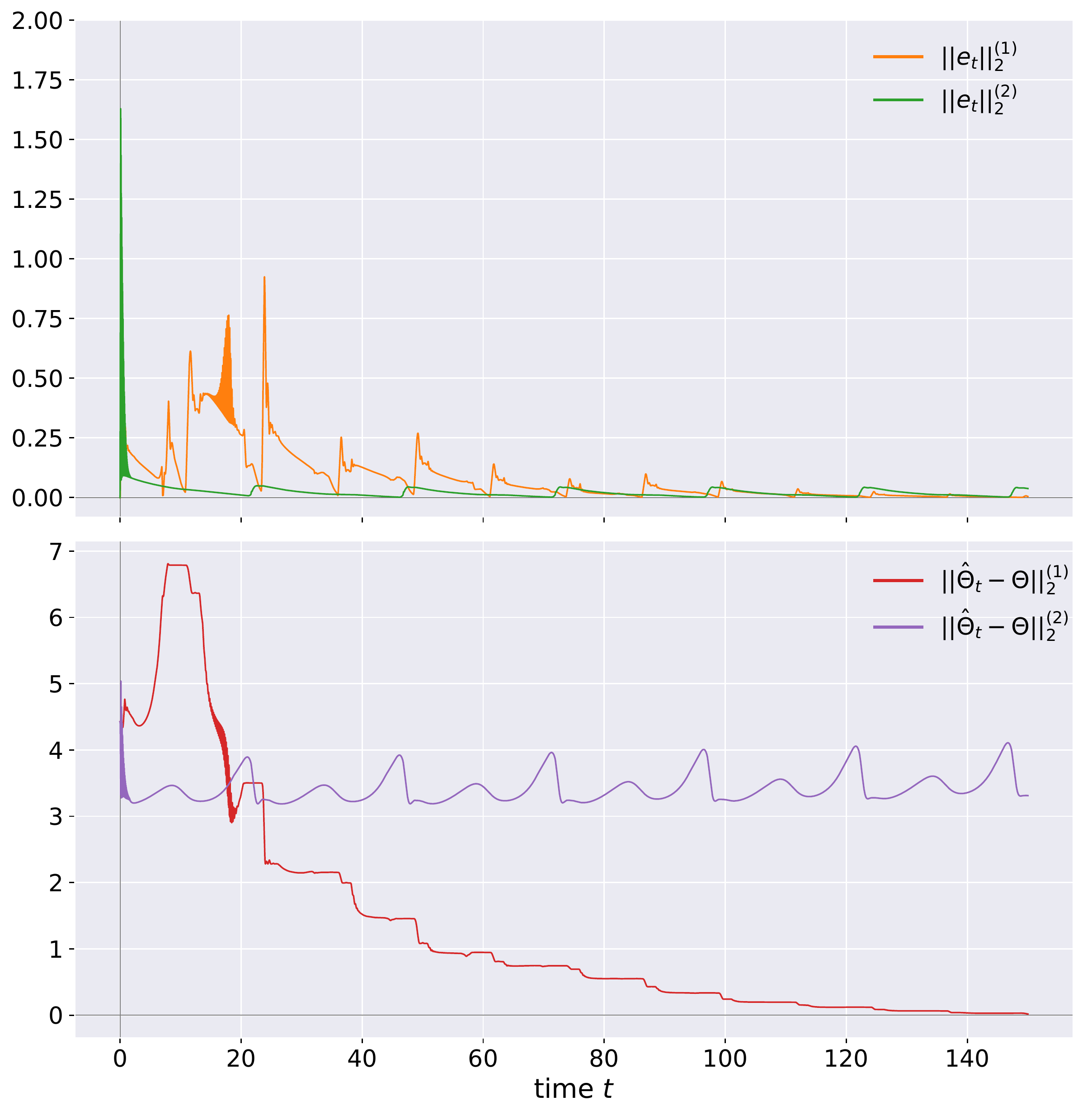}
\caption{\label{fig:thEstConv}Norm of state tracking error $\|e_t\|_2$ converges to zero in both
cases, while the parameter estimation error $\|\Thetah_t-\Theta\|_2$ converges to zero for the
case with $r_t^{(1)}$ and but not for the case with $r_t^{(2)}$.}
\end{figure}
\fi


\section{Conclusion}
\label{sec:conclusion}

In this paper, we defined a geometric criteria that leads to sufficient conditions for
persistency of excitation with single hidden-layer neural networks of step or ReLU activation
functions. Future work will focus on using the function approximation properties of ReLU
activations to obtain similar results when the plant nonlinearity is \textit{not known} but can
be \textit{approximated} by $\Theta^\top\phi(x_t)$ with ReLU's.

\if\ARXIV1 
\newpage
\phantom{ }
\newpage
\phantom{ }
\newpage
\fi


\printbibliography

\if\ARXIV1
\newpage
\appendices
\onecolumn
\section{Proofs of Theorems \ref{thm:step} and \ref{thm:ReLU}}\label{appThmProofs}
Here we will prove the main results provided in Section~\ref{sec:theory}.

\subsection{Proof Methodology and Setup}

The proofs for both theorems rely on the same contradiction method and share a common setup.
First, we will assume that the state trajectory $x_t$ is continuous and stays within some
compact set $\Bc\subset\R^n$ for all $t\geq0$, and that the component activations
$\sigma_1(w_1^\top x+b_1),\dots,\sigma_N(w_N^\top x+b_N)$ of $\phi$ are then bounded on that
$\Bc$.\\\\
Then we will assume that for any shift $\tau\geq0$ of the time window $[\tau,\tau+T]$ of any
length $T>0$, there exists a \underline{nonzero} vector $v\in\R^N\setminus\{0_N\}$ such that
\begin{equation}\label{app:nonPEvec0}
v^\top\left( \int_{\tau}^{\tau+T} \phi(x_t)\,\phi(x_t)^\top\dt\right) v \ = \ 0 \ \ .
\end{equation}
Finally, we will show that if the state trajectory $x_t$ meets certain requirements within
any shift $\tau\geq0$ of the window $[\tau,\tau+T^*]$, for some window length $T^*>0$, then in
fact \eqref{app:nonPEvec0} can only hold if $v=0_N$. This is a contradiction and thus proves that
\begin{equation}\label{strictPEPD}
\int_{\tau}^{\tau+T^*} \phi(x_t)\,\phi(x_t)^\top\dt 
\end{equation}
must be strictly positive definite, with bounded eigenvalues, uniformly over all windows
$[\tau,\tau+T^*]$. Then by the Min-max principle, the persistency of excitation conditions
\eqref{PE} and \eqref{PEvec} are satisfied with setting $\alpha_1$ and
$\alpha_2$ respectively as the smallest and largest eigenvalues of \eqref{strictPEPD} over
all shifts $\tau\geq0$ of the window $[\tau,\tau+T^*]$.\\\\
And so, let there be such a continuous trajectory $x_t$ staying within a compact set
$\Bc\subset\R^n$ for all $t\geq0$. Within any particular time window $[\tau,\tau+T]$, the
trajectory visits $\Ls$ activation regions $\Ac_{\As_1},\dots,\Ac_{\As_\Ls}$ and the
corresponding index sequences $\As_1,\dots,\As_\Ls$ and $\Is_1,\dots,\Is_{\Ls-1}$ are defined as
in Section~\ref{sec:setup}.\\\\
For any time window $[\tau,\tau+T]$ and over the corresponding sequence $s\in[\Ls]$ of visited
activation regions, recall the definition of the subsets $\Tc_s$ from (\ref{eq:timePartition}):
\begin{equation*}
\Tc_s = \{t\in[\tau,\tau+T] \ | \ x_t \in \Ac_{\As_s} \} \ .
\end{equation*}
Clearly the union of these subsets gives back the entire window as
\begin{equation*}
\bigcup_{s=1}^\Ls \Tc_s \ = \ [\tau,\tau+T] \ ,
\end{equation*}
meaning \eqref{app:nonPEvec0} is equivalent to
\begin{equation}\label{nonPEvec}
\int_{\tau}^{\tau+T} (v^\top\phi(x_t))^2\dt =
\sum_{s=1}^\Ls \int_{\Tc_s} (v^\top\phi(x_t))^2\dt = 0
\end{equation}
holding for all time windows $[\tau,\tau+T]$ and corresponding $\Tc_s$ subsets.\\\\
Since each integral in the sum is over $t\in\Tc_s$, and by definition we know that $x_t$ is only
within $\Ac_{\As_s}$ during that subset of the time window, then by definition only the
activation indices $i\in \Ss_{\As_s}$ of $\phi(x_t)$ will be nonzero, which gives
\begin{equation*}
v^\top\phi(x_t) = \sum_{i=1}^N v_i\,\phi(x_t)_i = \sum_{\Ss_{\As_s}} v_i\,\phi(x_t)_i
:= v_{\As_s}^\top\,\phi(x_t)_{\As_s} \ .
\end{equation*}
Here we denote $v_{\As_s}$ and $\phi(x_t)_{\As_s}$ respectively as the indices of $v$ and
$\phi(x_t)$ according to the active set $i\in \Ss_{\As_s}$, for each $s\in[\Ls]$.
This means \eqref{nonPEvec} is equivalent to
\begin{equation}\label{nonPEvec2}
\sum_{s=1}^\Ls \int_{\Tc_s} (v_{\As_s}^\top\phi(x_t)_{\As_s})^2\dt \ = \ 0
\end{equation}
holding for all time windows $[\tau,\tau+T]$ and corresponding $\Tc_s$ subsets. And since each of
these integrals is always nonnegative, \eqref{nonPEvec2} is equivalent to
\begin{equation}\label{nonPETtk}
v_{\As_s}^\top\phi(x_t)_{\As_s} = 0 \qquad \forall t\in\Tc_s, \forall s\in[\Ls]
\end{equation}
holding for all time windows $[\tau,\tau+T]$ and corresponding $\Tc_s$ subsets.

\newpage
\subsection{Proof of Theorem~\ref{thm:step}}
In this case, \eqref{nonPETtk} is equivalent to
\begin{equation}\label{nonPETtkStep}
v_{\As_s}^\top\phi(x_t)_{\As_s} = v_{\As_s}^\top c_{\As_s} = 0 \qquad \forall t\in\Tc_s, 
\forall s\in[\Ls]
\end{equation}
holding for all time windows $[\tau,\tau+T]$ and corresponding $\Tc_s$ subsets. Here we denote
$c_{\As_s}$ as the indices of $c$ according to the active set $i\in \Ss_{\As_s}$, for all
$s\in[\Ls]$. This gives a fixed, nonzero $|\Ss_{\As_s}|$ length vector.\\\\
Note we have assumed that $x_t$ over $t\geq0$ is such that there exists $T^*>0$ whereby all $N$
hyperplanes are crossed and only at nondegenerate borders, within any time window
$[\tau,\tau+T^*]$. Thus, we can combine \eqref{activeSetInduction} and
\eqref{nonPETtkStep} to get that
\begin{align}\label{nonPESkInductionStep}
&0 = v_{\As_{s+1}}^\top c_{\As_{s+1}} \qquad \forall
t\in\Tc_{s+1} \\
&= 
\begin{cases}
v_{\As_s}^\top c_{\As_s} + v_{\Is_s}c_{\Is_s} \qquad \text{if}\ \ 
\Ac_{\As_{s+1}}\subseteq X_{\Is_s}^+ \\
v_{\As_s}^\top c_{\As_s} - v_{\Is_s}c_{\Is_s}  \qquad \text{if}\ \ 
\Ac_{\As_{s+1}}\subseteq X_{\Is_s}^\circ
\end{cases}\nonumber
\end{align}
must hold over all $s\in[\Ls-1]$ with $v_{\As_1}^\top c_{\As_1}=0$ for all $t\in\Tc_1$, for all
time windows $[\tau,\tau+T^*]$ and corresponding $\Tc_s$ subsets. However, we will now show by
induction that this can only hold if $v=0_N$.\\\\
Starting at $s=1$: it must hold that $v_{\As_1}^\top c_{\As_1} = 0$ for all $t\in\Tc_1$, and
since $v$ and $c$ are fixed it must continue to hold for all $t\in\Tc_2$, during which it must
hold that $v_{\As_2}^\top c_{\As_2} = 0$. And so, these can both hold only if the additional term
is $v_{\Is_1}c_{\Is_1}=0$, and since $c_{\Is_1}\neq0$ it must be that $v_{\Is_1}=0$. We can
then inductively make the same argument over the remaining $s=2,3,\dots,\Ls-1$ to obtain
$v_{\Is_1}=\dots=v_{\Is_{\Ls-1}}=0$. And since we have assumed all $N$ hyperplanes were crossed
at nondegenerate borders over the $\Ls-1$ hyperplane crossings, we therefore must have that
$v_i=0$ for all $i\in[N]$. This establishes the contradiction from the original assumption of a
nonzero $v$. \\\\
Since we assume that all $N$ hyperplanes are crossed and only at nondegenerate borders
during any time window $[\tau,\tau+T^*]$, this contradiction must hold for all such windows.
Therefore, \eqref{strictPEPD} must be strictly positive definite for all $\tau\geq0$ with
this $T^*$ window length. And since the (scaled) step is bounded on any compact set
$\Bc\subset\R^n$, we have proven persistency of excitation holds for any window length of at
least $T^*$.
\vphantom{a}\hfill$\blacksquare$

\subsection{Proof of Theorem~\ref{thm:ReLU}}
In this case, \eqref{nonPETtk} is equivalent to
\begin{equation}\label{nonPETtkReLU}
v_{\As_s}^\top\phi(x_t)_{\As_s} = v_{\As_s}^\top (W_{\As_s}^\top x_t + b_{\As_s}) = 0 \qquad
\forall t\in\Tc_s, \forall s\in[\Ls]
\end{equation}
holding for all time windows $[\tau,\tau+T]$ and corresponding $\Tc_s$ subsets. Here we denote
$W^\top$ as the $N\times n$ matrix with $w_i^\top $ as the $i$th row and $b$ as the length $N$
vector with $b_i$ as the $i$th entry, for all $i\in[N]$, and thus the $\As_s$ subscripts are
the indices according to the active set $i\in \Ss_{\As_s}$, for each $s\in[\Ls]$.\\\\
Let $t_1,\hat t_1,\ldots,t_M,\hat t_M\in\Tc_s$ be the times from condition (iii). Then for
$i=1,\ldots,M$ it must hold by (\ref{nonPETtkReLU}) that
\begin{equation}\label{eq:differenceNull}
v_{\As_s}^\top (W_{\As_s}^\top x_{t_i}+ b_{\As_s}) -
v_{\As_s}^\top (W_{\As_s}^\top x_{\hat t_i}+ b_{\As_s}) = 0 =
v_{\As_s}^\top W_{\As_s}^\top (x_{t_i}-x_{\hat t_i}) \ .
\end{equation}
In particular, $v_{\As_s}$ must be in the nullspace of
$$
\begin{bmatrix}
  (x_{t_1}-x_{\hat t_1})^\top \\
  \vdots \\
  (x_{t_M}-x_{\hat t_M})^\top
\end{bmatrix}W_{\As_s}=:\hat W_{\As_s} \ .
$$
Condition (iii) implies that $\hat W_{\As_s}$ and $W_{\As_s}$ have the same rank, and so the
rank-nullity theorem implies that their nullspaces must have the same dimension. Now, since the
nullspace of $\hat W_{\As_s}$ is a subspace of the nullspace  of $W_{\As_s}$, the nullspaces
must actually be the same. As a result, (\ref{eq:differenceNull}) implies that $v_{\As_s}^\top
W_{\As_s}^\top =0$ and so $v_{\As_s}^\top b_{\As_s}=0$ as well.\\\\
We have just shown that with condition (iii) satisfied, (\ref{nonPETtkReLU}) is equivalent to
\begin{equation}\label{nonPETtkReLU2}
v_{\As_s}^\top [W_{\As_s}^\top\ b_{\As_s}] = 0_{n+1}^\top \qquad  
\forall t\in\Tc_s, \forall s\in[\Ls]
\end{equation}
holding for all time windows $[\tau,\tau+T]$ and corresponding $\Tc_s$ subsets.\\\\
Note we have assumed that $x_t$ over $t\geq0$ is such that there exists $T^*>0$ whereby all $N$
hyperplanes are crossed and only at nondegenerate borders, within any time window
$[\tau,\tau+T^*]$. Thus, we can combine \eqref{activeSetInduction} and
\eqref{nonPETtkReLU2} to get that
\begin{align}\label{nonPESkInductionReLU}
&0_{n+1}^\top = v_{\As_{s+1}}^\top [W_{\As_{s+1}}^\top\ b_{\As_{s+1}}]  \qquad 
\forall t\in\Tc_{s+1} \\
&= 
\begin{cases}
v_{\As_s}^\top [W_{\As_s}^\top\ b_{\As_s}] + v_{\Is_s}[w_{\Is_s}^\top\ b_{\Is_s}] \qquad
\text{if}\ \ \Ac_{\As_{s+1}}\subseteq X_{\Is_s}^+ \\
v_{\As_s}^\top [W_{\As_s}^\top\ b_{\As_s}] - v_{\Is_s}[w_{\Is_s}^\top\ b_{\Is_s}] \qquad
\text{if}\ \ \Ac_{\As_{s+1}}\subseteq X_{\Is_s}^\circ
\end{cases}\nonumber
\end{align}
must hold over all $s\in[\Ls-1]$ with $v_{\As_1}^\top [W_{\As_1}^\top\ b_{\As_1}]=0_{n+1}^\top$
for all $t\in\Tc_1$, for all time windows $[\tau,\tau+T^*]$ and corresponding $\Tc_s$
subsets. However, we will now show by induction that this can only hold if $v=0_N$.\\\\
Starting at $s=1$: it must hold that $v_{\As_1}^\top [W_{\As_1}^\top\ b_{\As_1}] = 0_{n+1}^\top$
for all $t\in\Tc_1$, and since $v$, $W$, and $b$ are fixed it must continue to hold for all
$t\in\Tc_2$, during which it must hold that $v_2^\top [W_{\As_2}^\top\ b_{\As_2}]=0_{n+1}^\top$.
And so, these can both only hold if the additional term is $v_{\Is_1}[w_{\Is_1}^\top\ b_{\Is_1}]
=0_{n+1}^\top$, and since $[w_{\Is_1}^\top\ b_{\Is_1}]$ is nonzero it must be that
$v_{\Is_1}=0$. We can then inductively make the same argument over the remaining
$s=2,3,\dots,\Ls-1$ to obtain $v_{\Is_1}=\dots=v_{\Is_{\Ls-1}}=0$. And since we have assumed all
$N$ hyperplanes were crossed at nondegenerate borders over the $\Ls-1$ hyperplane crossings, we
therefore must have that $v_i=0$ for all $i\in[N]$. This establishes the contradiction from the
original assumption of a nonzero $v$.\\\\
Since we assume that all $N$ hyperplanes are crossed and only at nondegenerate borders
during any time window $[\tau,\tau+T^*]$, this contradiction must hold for all such windows.
Therefore, \eqref{strictPEPD} must be strictly positive definite for all $\tau\geq0$ with
this $T^*$ window length. And since the ReLU is bounded on any compact set $\Bc\subset\R^n$, we
have proven persistency of excitation holds for any window length of at least $T^*$.
\vphantom{a}\hfill$\blacksquare$
\vphantom{a}\\\\

\subsection{Extending to Other Nonlinear, Positive Semidefinite Activations}

The reasoning and methods described in this appendix section should also be able to be applied
similarly to \textit{any} nonlinear, positive semidefinite activation function for which it is
possible to show the following:
\begin{enumerate}
\item describe all specific trajectories that $x_t$ can take such that
$\phi(x_t)_{\As_s}$ is some constant vector for the entire time within the activation region
$\Ac_{\As_s}$, and
\item formulate conditions like \eqref{nonPESkInductionStep} from
\eqref{nonPETtkStep} and \eqref{nonPESkInductionReLU} from
\eqref{nonPETtkReLU}, such that the condition is independent of the
specific value of $x_t$, and instead relies only on \textit{fixed} values within the definition
of the activation (like $c$ and $W$, $b$).\\
\end{enumerate}
If these two requirements are possible to meet, then the induction proof method will hold as
it does in this appendix section for (scaled) steps and ReLU activation functions.

\newpage
\section{Global Uniform Asymptotic Stability of \eqref{paramEstDyn}}\label{appLTVanderson}
Here we will show how the linear time-varying (LTV) system
\begin{equation}\label{appPELTV}
\dot{y}_t = \Gamma\,\phi(x_t)\phi(x_t)^\top y_t
\end{equation}
is proven globally uniformly asymptotically stable if the vector function $\phi:\R^n\to\R^N$
satisfies the persistency of excitation condition \eqref{PE}. State trajectories
$y:[0,\infty)\to\R^N$ and $x:[0,\infty)\to\R^n$ are assumed to be continuous (in time), and
$x_t$ is assumed to stay within some compact set $\Bc\in\R^n$ for all $t\geq0$. The vector
function is defined $\phi(x):=[g_1(x)\ \cdots\ g_N(x)]^\top$, where the piecewise continuous,
nonlinear scalar functions $g_1,\dots,g_N:\R^n\to\R$ are bounded on $\Bc$, and
$\Gamma$ is some symmetric, positive-definite $N\times N$ matrix with $\|\Gamma\|_F=\Gs$.
Then, $\phi(x_t)$ is a piecewise continuous (and regulated) vector function of time, mapping
$[0,\infty)\to\R^N$.\\\\
Firstly, Lemma 1 of \cite{anderson1977exponential} says that if $F:[0,\infty)\to\R^{N\times N}$,
$H:[0,\infty)\to\R^{N}$, and $K[0,\infty)\to\R^{N}$ are regulated vector or matrix functions of
time, and there exists positive scalars $\alpha_3,T>0$ such that
$$
\int_{\tau}^{\tau+T} \|K_t\|_2^2 \,\dt \ \leq\ \alpha_3 
$$
holds for all $\tau\geq0$, then $[F_t,H_t]$ is uniformly completely observable (UCO) if and only
if $[F_t+K_t\,H_t^\top,H_t]$ is UCO.\\\\
Next, Lemma 2 of \cite{anderson1977exponential} then says that the global uniform asymptotic
stability of
$$
\dot{y}_t = F_t\, y_t
$$
is equivalent to: i) the existence of a symmetric, differentiable matrix function of time
$P:[0,\infty)\to\R^{n\times n}$ satisfying, for some positive scalars $\beta_1,\beta_2>0$, that
\begin{equation*}
\beta_1I_N \preceq P_t \preceq \beta_2 I_N \qquad\text{and}\qquad
-\dot{P}_t = P_t\,F_t + F_t^\top P_t + H_tH_t^\top
\end{equation*}
both hold for all $t\geq0$, \textbf{and} ii) that $[F_t,H_t]$ is UCO. \\\\
Finally, a slight modification of Theorem 1 of \cite{anderson1977exponential}, in order to
allow for $\Gamma$ matrices other than the identity $I_N$, says that the global uniform
asymptotic stability of \eqref{appPELTV} holds if $\phi(x_t)$ satisfies the persistency of
excitation condition \eqref{PE}.\\[8pt]
\noindent\textit{Proof.} Suppose that $\phi(x_t)$ satisfies the persistency of excitation
condition \eqref{PE} for some positive scalars $\alpha_1,\alpha_2,T>0$. This means, by
definition, that $[0_{N\times N},\phi(x_t)]$ is UCO. Now set $K_t=-\Gamma\,\phi(x_t)$ and note
that \eqref{PE} gives
$$
\int_{\tau}^{\tau+T}\|K_t\|_2^2 \,\dt \ \leq\ 
\int_{\tau}^{\tau+T}\|\Gamma\|_F^2\|\phi(x_t)\|_F^2 \,\dt =
\Gs^2\int_{\tau}^{\tau+T}\! \tr\left[\phi(x_t)\phi(x_t)^\top\right] \dt \ \leq\ 
\Gs^2\alpha_2N \ ,
$$
thus satisfying Lemma 1 and so $[-\Gamma\,\phi(x_t)\phi(x_t)^\top,\phi(x_t)]$ is UCO.\\\\
Therefore, we have that Lemma 2 is satisfied with $P_t=\frac{1}{2}\Gamma^{-1}$,
$F_t=-\Gamma\,\phi(x_t)\phi(x_t)^\top$, and $H_t=\phi(x_t)$. This proves \eqref{appPELTV}
is globally uniformly asymptotically stable.
\vphantom{a}\hfill$\blacksquare$
\\
\begin{remark}
Lemma 1 and Theorem 1 of \cite{anderson1977exponential} use \textit{exponential} stability in
their statements, rather than global asymptotic uniform stability. However, for linear systems
these types of stability are equivalent. See Section 1.5.2 of \cite{sastry1989adaptive} for
details.
\end{remark}
\vphantom{a}
\begin{remark}
Lemma 1 and 2 and Theorem 1 of \cite{anderson1977exponential} also hold with $\phi(x_t)$,
$K_t$, and $H_t$ being regulated $N\times\ell$ matrix functions of time for any $\ell\geq1$,
using the induced Euclidean matrix norm when $\ell>1$. We specified their results for $\ell=1$
and the Euclidean vector norm, in order to match the application in
Section~\ref{sec:numerical}.
\end{remark}


\fi

\end{document}